\documentclass[10pt,titlepage]{article}
\usepackage{latexsym,amssymb,graphicx}
\usepackage{srcltx} 
\usepackage{latexsym} 
\usepackage{amsfonts} 
\usepackage{epsfig,float} 
 \restylefloat{figure}
\date{1999}

 \newtheorem{thm}{Theorem}[section]

\newcommand{\ds}{\displaystyle}
\def\be{\begin{equation} }
\def\ee{\end{equation} }

\begin{document}
\vspace{1cm}

\begin{center}

{\large {\bf A Brief Overview of the Sock Matching Problem}}

\renewcommand{\thefootnote}{}


\renewcommand{\thefootnote}%

\vspace*{5mm}

   Bojana Panti\'{c} and Olga Bodro\v{z}a-Panti\'{c}  \vspace*{3mm}

 Department of Mathematics and Informatics \\
  Faculty of Sciences  \\
  University of Novi Sad \\ Serbia

\vspace*{3mm}

 e-mail: dmi.bojana.pantic@student.pmf.uns.ac.rs \\
 olga.bodroza-pantic@dmi.uns.ac.rs

\vspace*{1cm}

\end{center} \vspace*{3mm}





\parbox{10cm}{{\bf Abstract.}   This short note deals with the so-called {\em
Sock Matching Problem}. We define $B_{n,k}$ as the number of all
the finite sequences $a_1, \ldots, a_{2n}$ of nonnegative integers
which contain at least one occurrence of $k$ $(1 \leq k \leq n)$
and for which $a_1 = 1 $,  $a_{2n}=0$ and $ \mid a_i -a_{i+1}\mid
\; = 1$. The value $a_i$ can be interpreted as the number of
unmatched socks being present after having drawn the first $i$
socks randomly out of the pile which initially contained $n$ pairs
of socks. Here, establishing a link between this problem and with
both some old and some new results, related to the number of
restricted Dyck paths, we obtain a few valid forms of the sock
matching theorem and  prove that the probability for $k$ unmatched
socks to appear (in the very process of drawing one sock at a
time) approaches $1$ as the number of socks becomes large enough.}

\vspace*{10mm}


     \section{Introduction}
In simple terms, what is understood under The Sock Matching
Problem \cite{GJRW} is the following procedure. Out of the laundry
pile that contains exactly n different pairs of socks socks are
being drawn randomly, one at a time (so that in the end all the
$2n$ socks get matched). In each move one tries to find the
adequate pair among the drawn socks, in case it had already been
obtained in the process. Furthermore, each of the two options:
drawing a match for some sock or drawing a sock that has no match
as of yet matches a single move, either one unit up or one unit to
the right,
 on a Dyck path (a path in an $n \times n$ grid starting from the lower left
corner $(0,0)$ and ending in the upper right corner $(n,n)$ using
merely moves up and to the right without ever crossing the
diagonal, see Fig. 1a) - this particular model of the Dick paths
was used in \cite{GJRW}).

It is a well known fact indeed that the number of all the Dick
paths of order $n$ is equivalent to the $n^{th}$ Catalan number
$\ds C_n = \frac{1}{n+1} {2n \choose n}$ (see \cite{S1}).

Let us now formulate our Sock Matching Problem in somewhat
mathematically stricter terms. In fact, let us focus upon the
total number of ways, labeled by $B_{n,k}$, to get at least $k$
unmatched socks at least once during the matching process.
Considering the aforementioned interpretation of our problem that relies on
the use of Dick paths, it is our task to determine which ones out of these
$C_n$ possibilities are those that present the paths which hit or
pass above the line $y=x+k$.

Admittedly there is a wide range of interpretations of the Catalan
numbers $C_n$. However, it is that which allows various authors to
opt themselves for the most suitable one. Here, we  make use of
the following terminology from \cite{S1,S2}:

\begin{itemize}
\item[-]
Lattice paths, used in \cite{II}, consider the up- and down-steps.
The former $(1,1)$-steps represent the case when "a sock with no
match has been drawn", whereas the latter $(1,-1)$-steps represent
the case when "a match has been made". Here, a lattice path goes
from $(0,0)$ to $(2n,0)$ on the Cartesian plane without ever
moving across the $x$-axis (though it is allowed to hit it), as
shown in Fig. 1b);

\item[-]
The number of planted plane trees (i.e. rooted trees which have
been embedded in the plane so that the relative order of subtrees
at each branch is part of its structure; ordered trees) with $n+1$
nodes, see Fig. 2;

\item[-]
Discrete random walks in a straight line with an absorbing barrier
at  $0$ (represented by the sequences $1=c_1, c_2,c_3, \ldots,
c_{2n+2}=0$; \ where  $c_i \geq 1$ for $i<2n+2$ and  $\mid c_i
-c_{i+1} \mid \; = 1$);

\item[-]
The number of all the sequences $a_1, \ldots, a_{2n}$ of
nonnegative integers with $a_1 = 1 $,  $a_{2n}=0$ and $ \mid a_i
-a_{i+1} \mid \; = 1$ (Problem 6.19 ($u^5$) in \cite{S2}). (Note
that the value $a_i$  can be interpreted as the exact number of
unmatched socks which are present after drawing the $i^{th}$ sock.
Let us also note that further on this interpretation shall be refer to
as the one with the nonnegative sequences).
\end{itemize}


\begin{figure}
  \includegraphics[width=4.5in]{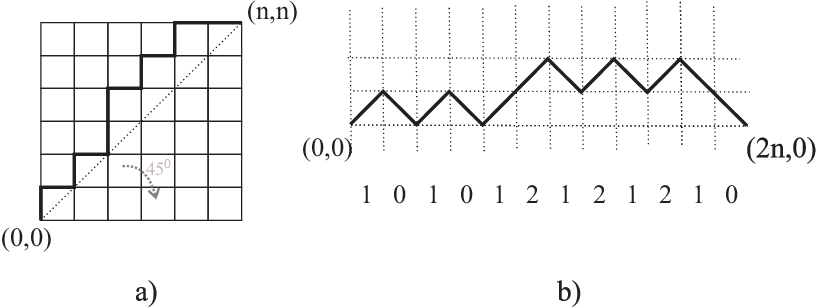}
\caption{a) a Dyck path;  \ b) a lattice path with a corresponding
nonnegative sequence}
\label{fig1}       
\end{figure}

A bijection is easily established between the sets of the
aforementioned combinatorial objects. For instance, the Dyck path
shown in Fig. 1a) corresponds to both the path in Fig. 1b) (which
can be obtained by rotating that very figure for $-45^0$ and then
expanding it with the expansion coefficient of $\sqrt{2}$) as well
as to the tree in Fig. 2. To be more precise, by wandering around
that tree the vertical component of successive positions describes
a path from $1$ (the root of the tree) to $0$. Consequently, in
this particular example the corresponding discrete random walk
would be $1,2,1,2,1,2,3,2,3,2,3,2,1,0$; whereas the nonnegative
sequence, mentioned in the last interpretation, would be
$1,0,1,0,1,2,1,2,1,2,1,0$.

The height of a Dyck path is the greatest distance from the
diagonal to the path, the height of a lattice path the greatest
distance from the $x$-axis to the path, whilst the height of a
planted ordered tree is the number of nodes on a maximal simple
path starting at a root. Clearly, the height of a Dyck path is
$\ds\max_i a_i$ in the nonnegative sequences interpretation,
whereas the height of a corresponding planted plane tree is
$\ds\max_i c_i$ in the discrete random walks interpretation. Now,
it is worthwhile realising that the latter value is for one
greater than the former. To illustrate this point, take another
look at Figures 1 and 2. There, the height of the presented Dyck
path (Fig. 1) is 2 (at most 2 unmatched socks appear), as opposed
to the height of the planted plane tree (Fig. 2) which is 3!

Bearing all this in mind, the outlined problem from the heading,
i.e. the number $B_{n,k}$ may represent as follows:
\begin{itemize}
\item[-] the number of Dyck paths of height at least $k$ (the ones that hit or cross the line $y=x+k$),
\item[-] the number of lattice paths of height at least  $k$ (the ones that hit or cross the line $y=k$),
\item[-] the number of planted plane trees of height at least $k+1$,
\item[-] the number of discrete random walks with  $\ds \max_i c_i \geq k+1$,
\item[-] the number of all the nonnegative  sequences containing the letter $k$.
\end{itemize}


\begin{figure}[htbp]
\begin{center}
\includegraphics[width=2in]{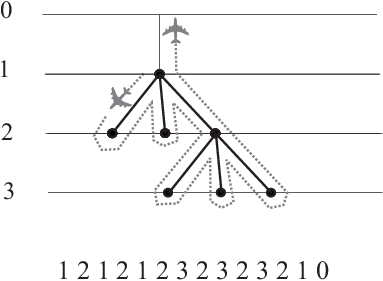}
\end{center}
\caption{a tree with its random walk} \label{fig2}
\end{figure}

Contemporary researches related to the Dyck paths refer mainly to
the number of restricted lattice paths, where crossing the
$x$-axis is allowed. Ili\' c and Ili\'c in \cite{II} gave the upper
and lower bounds for this number in the form of binomial
coefficients. Forging a link between this problem and an older
paper \cite{PP}
 from 1985. H. Prodinger in \cite{P} provides an explicit formula for those, as given in this theorem:

\begin{thm}
The number of random walks from $(0,0)$ to $(2n,0)$ with up-steps
and down-steps of one unit each, under the condition that the path
is placed between the lines $y=-h$ and $y=t$ is equal to
\\
$ \ds  \sum_{ j \geq 0} \left[  {2n  \choose n-j(h+t+2)} - {2n
\choose n-j(h+t+2)-h-1}  \right.$ \be \label{0}  - {2n \choose
n-j(h+t+2)-t-1} + \left.
 {2n  \choose n- (j+1)(h+t+2)}  \right] .\ee
\end{thm}

 In the special case, for  $h=0$, we obtain the number of all the sequences $a_1, \ldots, a_{2n}$ of nonnegative integers  with
$a_1 = 1 $,  $a_{2n}=0$ and $ \mid a_i -a_{i+1} \mid \; = 1$ and
$a_i \leq t$ $(n \geq 1, t \geq 0)$, which we label by $A_{n,t}$.
Obviously, $A_{n,0}=0$ and $A_{n,1}=1$ for $n \geq 1$,
$A_{n,t}=C_n$ for $t \geq n$. The value $A_{n,t}$ was already
essentially obtained in \cite{BKR} in the distant 1972. in the
form of trigonometric functions. The authors of that paper used
the rooted tree (planted plane tree) interpretation. As for
convenience, we reformulate their results in the following
theorem.

\begin{thm}
\be \label{2} \ds A_{n,t} = \frac{1}{t+2} \sum_{1 \leq j \leq
\frac{t+1}{2}} 4^{n+1} \sin^2\left( \frac{j \pi}{t+2} \right)
\cos^{2n} \left( \frac{j \pi}{t+2} \right),  \ee where $ n \geq
1$ and $t \geq 1$.
 \end{thm}

It is quite an interesting fact that this formula has been
rediscovered many a time, as the authors of paper
\cite{BKR} clearly point out, and that above all Lagrange derived
a formula in 1775. which essentially includes this as a special
case.

The authors of  \cite{GJRW} derive a recurrence formula, to which
they refer to as the so-called {\em Sock Matching Theorem}, for
the numbers  $B_{n,k}$. Additionally, they make a proposition that
the probability $P_{n,k}$ for a Dyck path to reach the line
$y=x+k$ approaches $1$ as the number n becomes large enough, i.e.
$\ds \lim_{n \rightarrow \infty} \frac{B_{n,k}}{C_{n}} = 1$.
However, we noticed some inaccuracies in these proofs (for more
details on that see \cite{PB}, page 3 and 4). For the purpose of
providing a valid formula for $B_{n,k}$, with the accompanying
proofs, we shell present two equivalent alternatives in Section 2,
using more than just the mentioned authors' idea. In section 3 we
give an explicit expression for  $B_{n,k}$. In Section 4 we prove
the mentioned theorem about the asymptotic behavior of the ratio
$B_{n, k}/C_n$ when n converges to infinity.

\section{The sock matching theorem}

\begin{thm}{\bf (The Sock Matching Theorem - I alternative)}

The sequence $B_{n,k}$ whose $n^{th}$ term represents the number
of Dyck paths of order $n$ which hit or cross the line $y = x + k$
is determined by the following recurrence formula: \be \label{2.1}
B_{n,k} = \sum_{i=1}^{n} (B_{i-1,k-1}C_{n-i} + C_{i-1}B_{n-i,k} -
B_{i-1,k-1}B_{n-i,k}).\ee
\end{thm}

\noindent {\bf Proof (The first one):} Let $(i,i)$ be the first
point on the line $y=x$ which the Dyck path visits after $(0,0)$.
Further, we take three possibilities into consideration: the line
hits $y=x+k$ before $(i,i)$ (Case 1), the line hits $y=x+k$ after
$(i,i)$ (Case 2), and the line hits $y=x+k$ both before and after
$(i,i)$ (Case 3).

The number of paths in the first case is $B_{i-1,k-1}C_{n-i}$.
Namely, the number of ways to hit $y=x+k$ between  $(0,0)$ and
$(i,i)$ without hitting $y=x$ is the same as the number of ways to
get from $(0,1)$ to $(i-1,i)$ hitting  $y=x+k$ but not crossing
$y=x+1$, which is $B_{i-1,k-1}$ and the number of ways to get from
$(i,i)$ to $(n,n)$ without crossing $y=x$ is $C_{n-i}$.

Similarly, the numbers in the second and third case are
$C_{i-1}B_{n-i,k}$ and $B_{i-1,k-1}B_{n-i,k}$, respectively. $
\Box $


\noindent {\bf Proof (The second one):} There is, however, yet
another approach which may be made in order to obtain the formula
(\ref{2.1}) which makes use of the recurrence relation satisfied
by the numbers $A_{n,k}$ derived in \cite{BKR}: \be \label{2.3}
A_{n,k+1} = A_{n-1,k+1}A_{0,k} + A_{n-2,k+1}A_{1,k}+ \ldots +
A_{0,k+1}A_{n-1,k}, \ee where $ n \geq 1 $ and $  k \geq 0 $; with
the initial conditions for $A_{n,0}=0, \; \; \mbox{when } \; n
\geq 1 $
 and for  $A_{0,k} = \hspace*{-4mm} ^{\rm def} \; \; 1 , \; \; \mbox{when } k \geq 0 $.
To be more specific, as \be \label{2.4} B_{n,k} = C_n - A_{n,k-1}
, \ee  where $ n \geq 1 $ and $ k\geq 1 $ (evidently, $B_{n,1}
=C_n$) making the necessary substitutions in (\ref{2.3}) and with
the use of the well-known recurrence relation
for Catalan numbers \\
($ \ds  C_0=1, \; \; C_{n+1} = \sum_{i=0}^{n} C_i C_{n-i}, \mbox{
for } n \geq 0 $) we obtain the desired relation. $ \Box $

\begin{thm}[The Sock Matching Theorem - II alternative]
\be \label{2.2}  B_{n,k} = \sum_{j=0}^{n-1} (B_{j,k}C_{n-j-1} +
C_{j}B_{n-j-1,k-1} - B_{j,k}B_{n-j-1,k-1}).\ee
\end{thm}

\noindent {\bf Proof:}  Similarly to the previous alternative we
take the point $(i,i)$ into consideration, only this time as the
last point on the line $y=x$ that the Dyck path visits before
$(n,n)$. Seen from this perspective, the corresponding numbers for
the three cases would be exactly the values $B_{j,k}C_{n-j-1} $, $
C_{j}B_{n-j-1,k-1}$ and $ B_{j,k}B_{n-j-1,k-1}$.

By the way, the proof for the second formula (\ref{2.2}) could,
obviously, be obtained from (\ref{2.1}) by a fairly simple
substitution: $i=n-j$. $ \Box $

\section{The explicit formula for $B_{n,k}$}

We now give the explicit expression for the values of  $B_{n,k}$.

\begin{thm}
\be \label{2.6} B_{n,k}=  \sum_{ j = 1}^{\lfloor \frac{n + 1}{k +
1}\rfloor}   {2  n + 2 \choose  n + 1 - j (k + 1)}  -4 \sum_{ j =
1}^{\lfloor\frac{n }{k + 1}\rfloor}
 {2  n  \choose  n  - j (k + 1)}.
\ee
\end{thm}

\noindent {\bf Proof (The first one):} Recall that the value
$A_{n,t}$ for the lower bound $h=0$ is known from (\ref{0}).
Further, setting the upper bound to be $ t=k-1$ it follows
directly from (\ref{2.4}) that \\
$ \ds B_{n,k} = \frac{1}{n+1}{2n \choose n}   - \sum_{ j \geq 0}
\left[  {2n  \choose n-j(k+1)} -   {2n  \choose n-j(k+1)-1}
  \right.
  $
 \be \label{2.5}  \left.  -   {2n  \choose
n-j(k+1)-k} +
 {2n  \choose n- (j+1)(k+1)}  \right] .\ee

After some minor algebraic simplifications of the above expression
we have \\
$\ds B_{n,k} =$
 \be \label{2.11}   \sum_{ j \geq 1} \left[  {2n \choose (n+1)-j(k+1)}
-2  {2n  \choose n-j(k+1)}   + {2n \choose
(n-1)-j(k+1)} \right],\ee  \\
which coincides with the formula presented in \cite{BKR} for the
number of planted plane trees with $n+1$ nodes whose height is
greater than $k$. Now, since \\ \\$ \ds {2  n + 2 \choose  n + 1 -
j (k + 1)}= {2  n + 1 \choose  n - j (k + 1)}+ {2  n + 1 \choose n
+ 1 - j (k + 1)}$
\\  \\ \\ $ \ds
={2  n  \choose  n -1 - j (k + 1)}+ 2{2  n \choose  n - j (k + 1)}
+{2  n \choose  n +1 - j (k + 1)}, $ \\ \\
 the expression (\ref{2.6}) is easily derivable from (\ref{2.11}).
$ \Box $


\noindent {\bf Proof (The second one):} Let us now commence from
formula (\ref{2}). By applying $\sin^2 \alpha = 1-\cos^2 \alpha$
to it we obtain the following \be \label{2.8} \ds A_{n,t} =
\frac{4^{n+1}}{t+2} \left( \sum_{1 \leq j \leq \frac{t+1}{2}}
\cos^{2n} \left( \frac{j \pi}{t+2} \right) - \sum_{1 \leq j \leq
\frac{t+1}{2}}   \cos^{2n+2} \left( \frac{j \pi}{t+2} \right)
\right), \ee where  $   n \geq 1 .$  Further, using one of the
notations for the representation of trigonometric power sums,
namely the one with binomial coefficients, from \cite{FGK}, we
have \be \label{2.9}\sum_{j=0}^{N-1}
\cos^{2m}\left(\frac{j\pi}{N}\right) = 2^{1-2m} N \left({2m-1
\choose m-1} + \sum_{p=1}^{\lfloor \frac{m}{N} \rfloor}  {2m
\choose m-p N} \right),\ee   where $ m,N \in \mathbb{N}$ and $ m
\geq N$. Now, making the necessary substitutions, i.e. $N=t+2$ and
for $m$ at first $m=n$ and then $m=n+1$, a brief simplification
process leads to
\\
 $\ds A_{n,t} = 4 \left[ {2n -1 \choose n-1}  + \sum_{ j \geq 1}^{\lfloor \frac{n}{t+2} \rfloor}    {2n  \choose n-j(t+2)}
  \right]$
  \be \label{2.10} -   \left[ {2n +1 \choose n}  +  \sum_{ j \geq 1}^{\lfloor \frac{n+1}{t+2} \rfloor}    {2n +2 \choose n+1-j(t+2)}
 \right]. \ee \\
Once again, utilising (\ref{2.4}) and yet again simplifying the
obtained expression we eventually come to the desired formula
(\ref{2.6}).  $ \Box $

\scriptsize \vspace*{0.3cm}
\begin{tabular}{||c||r|r|r|r|r|r|r||} \hline\hline
$B_{n,k}$ & $k=1$ & $2$  & $3$ & $4$ & $5$ & $6$ & $7$  \\
\hline\hline $B_{1,k}$ &1 &  0 & 0 & 0 & 0 & 0 & 0  \\
\hline       $B_{2,k}$ &2 & 1 & 0 & 0 & 0 & 0 & 0  \\ \hline
             $B_{3,k}$ &5 & 4 & 1 & 0 & 0 & 0 & 0\\ \hline
             $B_{4,k}$ &14 & 13 & 6 & 1 & 0 & 0 & 0 \\ \hline
             $B_{5,k}$ &42 &41 & 26 & 8 & 1 & 0 & 0 \\ \hline
             $B_{6,k}$ &132& 131 & 100 & 43 & 10 & 1 & 0 \\
\hline       $B_{7,k}$ &429 & 428 & 365 & 196 & 64 & 12 & 1 \\
\hline $B_{8,k}$ &1430 & 1429 &
1302 & 820 & 336 & 89 & 14 \\
\hline $B_{9,k}$ &4862 & 4861 & 4606 & 3265 & 1581 & 528 & 118 \\
\hline $B_{10,k}$ &16796 & 16795 & 16284 & 12615 & 6954 & 2755 &
780 \\ \hline
  $B_{11,k}$ &58786 & 58785 & 57762 & 47840 & 29261 & 13244 & 4466 \\ \hline
   $B_{12,k}$ &208012 & 208011 & 205964 & 179355 & 119438 & 60214 &   23276 \\ \hline
   $B_{13,k}$ &742900 & 742899 & 738804 &   667875 & 477179 & 263121 & 113620 \\ \hline
   $B_{14,k}$ &2674440 & 2674439 & 2666248 & 2478022 & 1877278 & 1116791 & 528840 \\ \hline
   $B_{15,k}$ &9694845 & 9694844 &   9678461 & 9180616 & 7303360 & 4637476 & 2375101
\\ \hline \hline
\end{tabular}

\scriptsize \vspace*{0.3cm}
\begin{tabular}{||c||r|r|r|r|r|r|r|r||} \hline\hline
$B_{n,k}$ & $k=8$ & $9$ & $10$ & $11$ & $12$ & $13$ & $14$ & $15$
\\ \hline\hline $B_{1,k}$ & 0 & 0 & 0
& 0 & 0 & 0 & 0 & 0 \\ \hline $B_{2,k}$ & 0 & 0 & 0 & 0 & 0 & 0 &
0 & 0 \\ \hline $B_{3,k}$ & 0 & 0 & 0 & 0 & 0 & 0 & 0 & 0\\ \hline
$B_{4,k}$  & 0 & 0 & 0 & 0 & 0 & 0 & 0 & 0\\ \hline $B_{5,k}$  & 0
& 0 & 0 & 0 & 0 & 0 & 0 & 0\\ \hline $B_{6,k}$  & 0 & 0 & 0 & 0 &
0 & 0 & 0 & 0\\ \hline $B_{7,k}$  & 0 & 0 & 0 & 0 & 0 & 0 & 0 & 0\\
\hline $B_{8,k}$  & 1 & 0 & 0 & 0 & 0 & 0 & 0 & 0\\
\hline $B_{9,k}$  & 16 & 1 & 0 & 0 & 0 & 0 & 0 & 0\\ \hline
$B_{10,k}$ & 151 & 18 & 1 & 0 & 0 & 0 & 0 &    0\\ \hline
  $B_{11,k}$ & 1100 & 188 & 20 &    1 & 0 & 0 & 0 & 0\\ \hline
   $B_{12,k}$  & 6854 & 1496 & 229 & 22 & 1 & 0 & 0 & 0\\ \hline
   $B_{13,k}$  & 38480 & 10075 & 1976 & 274 & 24 & 1 & 0 &   0\\ \hline
   $B_{14,k}$  &   200655 & 60606 & 14301 & 2548 & 323 & 26 & 1 & 0\\ \hline
   $B_{15,k}$  & 990756 & 336168 & 91756 &   19720 & 3220 & 376 & 28 & 1
\\ \hline \hline
\end{tabular}
\normalsize

Tabular 1: Numerical values of $B(n,k)$ for small values of $n$
and $k$

\section{Asymptotic behavior}

\begin{thm}
The probability $P_{n,k}$ of reaching a given fixed $k$ ($k \geq 1$) approaches
$1$ as $n$ approaches infinity, i.e.
 $$   \lim_{n \rightarrow \infty} P_{n,k} =  \lim_{n \rightarrow \infty} \frac{B_{n,k}}{C_{n}} = 1 .           $$
 \end{thm}

\noindent {\bf Proof:} Firstly, notice that the theorem trivially
holds for $k=1$. Thus, we shall assume that $k\geq 2$. Exploiting
(\ref{2.4}) and (\ref{2}) some more we have
$$\lim_{n \rightarrow \infty} P_{n,k} =  \lim_{n \rightarrow \infty} \frac{B_{n,k}}{C_{n}} = 1- \lim_{n \rightarrow \infty}
\frac{A_{n,k-1}}{C_{n}}.    $$
 $$ = 1 -  \lim_{n \rightarrow \infty}  \frac{ \ds \frac{1}{k+1}  \ds
  \sum_{1 \leq j \leq \frac{k}{2}} 4^{n+1} \sin^2\left( \frac{j \pi}{k+1} \right)  \cos^{2n} \left( \frac{j \pi}{k+1} \right)}{ \ds \frac{1}{n+1}  \frac{(2n)!}{n!
  n!}}.           $$
Since $ \ds 0 < \sin^2\left( \frac{j \pi}{k+1} \right) <1$ and
$\ds 0 < \cos \left( \frac{j \pi}{k+1} \right) \leq \cos \left(
\frac{\pi}{k+1} \right)$ for every $j$ such that $ 1 \leq j \leq
\ds \frac{k}{2}$ we have
$$ 1  \geq \lim_{n \rightarrow \infty} P_{n,k} \geq  1 -  \lim_{n \rightarrow \infty}  \frac{ \ds \frac{4^{n+1}}{k+1}  \cdot   \frac{k}{2}  \cos^{2n}
\left( \frac{\pi}{k+1} \right)}{ \ds \frac{1}{n+1}
\frac{(2n)!}{n! n!}} .$$
 Applying the Stirling's approximation we obtain the following
$$ \lim_{n \rightarrow \infty}  \frac{ \ds \frac{4^{n+1}}{k+1}  \cdot   \frac{k}{2}  \cos^{2n}
\left( \frac{\pi}{k+1} \right)}{ \ds \frac{1}{n+1}  \frac{(2n)!}{n! n!}} = \lim_{n \rightarrow \infty}
 \frac{ \ds 2k(n+1)\sqrt{\pi n} \cos^{2n}\left( \frac{\pi}{k+1} \right) }{k+1} . $$
Bearing in mind that $\ds \cos^{2n} \left( \frac{\pi}{k+1}
\right)$ approaches zero faster than $\ds
\frac{1}{(n+1)\sqrt{n}}$, it follows immediately that the last
limit is zero, leaving $\ds \lim_{n \rightarrow \infty} P_{n,k} =
1 $.
 $ \Box $


\subsection*{Acknowledgment}
We wish to express our sincerest gratitude towards Dragan
Stevanovi\'c for pointing out references \cite{P} and \cite{FGK},
as well as towards the anonymous referee for the thoughtful and constructive
remarks which helped improve the text.

Research supported by Grants OI 174018 and III 46005 of
the Ministry of Education and Science of
the Republic of Serbia



\begin{thebibliography}{1}
\bibitem{BKR} N. G. de Bruijn, D. E. Knuth, S. O. Rice, \textit{ The average height of planted plane trees, in
 Graph theory and computing}, edited by R.C.Read, Academic Press, New York
 15-22(1972)

\bibitem{GJRW}
 S. Gilliand, C. Johnson, S. Rush, D.  Wood, \textit{ The sock
matching problem},    Involve, \textbf{ 7}(5) (2014), 691--697.

\bibitem{FGK}C. M. da Fonseca, M. L. Glasser, V. Kowalenko, \textit{ Basic trigonometric
power sums with applications}, Ramanujan J. 42 (2017), 401--428.

\bibitem{II}
A. Ili\' c, A. Ili\' c, \textit{ On the number of restricted Dyck
paths},
 Filomat, \textbf{ 25}(3) (2011), 119--201.

\bibitem{PP}
W. Panny, H. Prodinger, \textit{  The expected height of paths for
several notions of height},  Studia Scientiarum Mathematicarum
Hungarica \textbf{ 20} (1985), 119--132.

\bibitem{PB}B. Panti\' c, O. Bodro\v za-Panti\' c, \textit{A Brief Overview of the Sock Matching Problem},  arXiv:1609.08353v1 [math.CO] 27 Sep 2016


\bibitem{P}
 H. Prodinger, \textit{The number of restricted lattice paths revisited},  Filomat, \textbf{ 26} (6)  (2012), 1133--1134.

\bibitem{S1}
R. P. Stanley, \textit{Enumerative Combinatorics}, Vol. I,
Cambridge University Press, Cambridge, 2002.

\bibitem{S2}
R. P. Stanley, \textit{Catalan Addendum to Enumerative Combinatorics},
 Volume 2,  version of 25 May 2013,
http://www-math.mit.edu/$\sim$rstan/ec/catadd.pdf, 2013.
\end{thebibliography}
     \end{document}